\documentclass[11pt,letterpaper]{amsart}

\usepackage{tikz}
\tikzstyle{every node}=[circle, draw, fill=black!50,
                        inner sep=0pt, minimum width=3pt]

\usepackage{amsfonts,amsmath,latexsym,color,epsfig,hyperref,enumitem, amssymb,bbm}

\newtheorem{theorem}{Theorem}[section]
\newtheorem{lemma}{Lemma}[section]

\renewcommand{\ge}{\geqslant}
\renewcommand{\leq}{\leqslant}
\renewcommand{\geq}{\geqslant}

\def\qed{\ifvmode\mbox{ }\else\unskip\fi\hskip 1em plus 10fill$\Box$}

\input{epsf}

\makeatletter
\def\Ddots{\mathinner{\mkern1mu\raise\p@
\vbox{\kern7\p@\hbox{.}}\mkern2mu
\raise4\p@\hbox{.}\mkern2mu\raise7\p@\hbox{.}\mkern1mu}}
\makeatother

\def\F{\mathbb F}

\def\F{\mathbb{F}}

\title{On skew corner-free sets}

\author{Cosmin Pohoata}
\thanks{Department of Mathematics, Emory University, Atlanta, GA. Email: {\tt cosmin.pohoata@emory.edu}. Research supported by NSF Award DMS-2246659.}

\author{Dmitrii Zakharov}
\thanks{Department of Mathematics, Massachusetts Institute of Technology, Cambridge, MA. Email: {\tt  zakhdm@mit.edu}. Research supported by the Jane Street Graduate Fellowship.}

\date{}

\begin{document}
\maketitle

\begin{abstract}
We construct skew corner-free sets in $[n]^2$ of size $n^{5/4}$, thereby disproving a conjecture of Kevin Pratt. We also show that any skew corner-free set in $\mathbb{F}_{q}^{n} \times \mathbb{F}_{q}^{n}$ must have size at most $q^{(2-c)n}$, for some positive constant $c$ which depends on $q$. 
\end{abstract}


\section{Introduction}

Motivated by matrix multiplication algorithms, Pratt \cite{pratt} asked the following nice question: what is the largest subset of $[n]^2$ which does not contain `skew corners', i.e. triples of points of the form
\begin{equation}\label{skew}
(x, y), (x, y+d), (x+d, y')    
\end{equation}
for $d \neq 0$ and $x, y, y'$ arbitrary?

We call such sets {\it{skew corner-free sets}}. Putting $y=y'$ in the above, we get that such a set does not contain {\it{regular}} corners $(x, y), (x+d, y), (x, y+d)$. Determining the maximum size of a subset of $[n]^2$ without regular corners is a well-studied problem in additive combinatorics. See for example \cite{green} and the references therein for some background. Given this connection to regular corner-free sets, Shkredov's result from \cite{shkredov} immediately implies that a skew corner-free set $S \subset [n]^2$ satisfies $|S| = O\left(n^2 / (\log \log n)^{c}\right)$, for some absolute constant $c > 0$.  On the other hand, Petrov \cite{petrov} constructed a skew corner free set in $[n]^2$ of size $\Omega( n \log n / \sqrt{\log \log n})$, and Pratt \cite{pratt} conjectured that any such set in $[n]^2$ has size $O(n^{1+\varepsilon})$ for any $\varepsilon > 0$. Furthermore, Pratt showed that such a result would show that certain approaches to matrix multiplication cannot achieve a running time of $O(n^{2+\varepsilon})$. We provide a construction which disproves this prediction.

\begin{theorem} \label{constr}
    There exists a skew corner-free set $S \subset [n]^2$ of size $\Omega(n^{5/4})$.
\end{theorem}

The main idea is to take advantage of a well-known property of the (affine version) of the classical Hermitian unital over $\mathbb{F}_{p^2}^{2}$, object which was also used by Mattheus--Verstraete \cite{mattheus2023asymptotics} in the recent breakthrough lower bound construction for the Ramsey number $R(4,t)$. We discuss the proof of Theorem \ref{constr} in Section 2. 

While Theorem \ref{constr} shows skew corner-free sets in $[n]^2$ may not always have size $O(n^{1+\varepsilon})$ for any $\varepsilon > 0$, we believe that such sets should in the very least still have size $O(n^{2-c})$ for some absolute constant $c > 0$. Establishing this fact already seems like an interesting problem, as there exist standard examples of corner-free sets in $[n]^2$ of size $n^{2-o(1)}$. However, it does not seem that the usual Fourier analytic methods can take advantage of the stronger condition of the set being skew corner free in any significant manner. In Section 3, we show a result in this spirit for the finite field model of this problem. 
\begin{theorem} \label{skewCLP}
    Let $q \ge 2$ be a prime and let $S \subset \mathbb{F}_{q}^{n} \times \mathbb{F}_{q}^{n}$ be any set without triples of the form $(\ref{skew})$
    with $x,y,y' \in \mathbb{F}_{q}^{n}$ and $d \in \mathbb{F}_{q}^{n} \setminus {0}$. Then, 
    $$|S| \leq 3q^{(2-c_{q})n},$$
    where the exponent $c_{q}$ is a positive constant depending on $q$ defined as 
    $$q^{1-c_{q}} = \inf_{0< x < 1} x^{-(q-1)/3}(1+x+\ldots+x^{q-1}).$$
\end{theorem}
Here we think of $q$ as fixed and as $n$ going to infinity. As the definition of $c_{q}$ might already suggest to the experienced reader, the proof of Theorem \ref{skewCLP} will use the so-called Croot-Lev-Pach lemma, famously introduced in \cite{clp17}, together with some of the ideas of Ellenberg and Gijswijt from their resolution of the cap set problem \cite{eg}. Qualitatively speaking, it is perhaps important to highlight that Theorem \ref{skewCLP} serves as a certificate that the skew corner-free problem in $\mathbb{F}_{q}^{n} \times \mathbb{F}_{q}^{n}$ does not obey the induced matching barrier described in \cite{barrier}, where Christandl, Fawzi, Ta, and Zuiddam show that the recent polynomial method as long as related tensor methods for upper bounding the Shannon capacity
(including slice rank, subrank, analytic rank, geometric rank, and G-stable rank) cannot yield a similar bound for the regular corner-free sets in $\mathbb{F}_{q}^{n} \times \mathbb{F}_{q}^{n}$.

\section{Proof of Theorem 1.1}
 Let $p \sim n^{1/4}$ be a prime and let $q = p^2$. For $a \in \F_q$ let $\bar a= a^p$ be the Galois conjugate and let $N(a) = a \bar a = a^{p+1}$ be the norm. Consider the affine version of the {\em Hermitian unital} in $\F_q^2$:
    $$
    Q = \{ (a, b) \in \F_q^2:~ N(a) + N(b) = 1 \}.
    $$
    It is well-known that $|Q| \sim q^{3/2}$ and for each point $x \in Q$ there exists a `tangent' $\F_q$-line $\ell_x \subset \F_q^2$ such that $\ell_x \cap Q = \{x\}$. In other words, $Q$ forms a so-called {\it{Nikodym set}} in $\mathbb{F}_{q}^{2}$. Indeed, if $x = (a, b)$ lies on $Q$ then we can define
    $$
    \ell_x = \{(a + t \bar b, b - t \bar a), ~ t \in \F_q\},
    $$
    where $\bar a = a^p$ denotes the Galois conjugate. Then we can write
    $$
    N(a + t \bar b) + N(b - t \bar a) = (a+t \bar b)(\bar a + \bar t b) + (b - t \bar a)(\bar b - \bar t a) = (1+N(t))(N(a)+N(b)) = 1+N(t)
    $$
    and so the point $(a + t \bar b, b - t \bar a)$ belongs to $Q$ if and only if $N(t) = 0$. But $N(t) = t^{p+1} = 0$ implies $t = 0$ finishing the proof of the claim. 
    
    Now we put
    $$
    S' = \{(x, y) \in \F_q^2 \times \F_q^2:~ x \in Q, ~ y \in \ell_x\}.
    $$
    Note that $|S'| = |Q| q \sim q^{5/2}$. Suppose that $S'$ contains a triple of the form $(\ref{skew})$, i.e. we have points $x, x+d \in Q$, $y, y+d \in \ell_x$ and $y' \in \ell_{x+d}$ for some $d \neq 0$. Then note that $x+d = x + (y+d) - y$, i.e. $x+d$ is an affine combination of points lying on $\ell_x$ and so $x+d$ itself has to lie on $\ell_x$. But then this contradicts the property that $Q \cap \ell_x = \{x\}$. So $S'$ is a skew corner-free set in $\F_q^2 \times \F_q^2$. 
    
    Now let us view $\F_q^2$ as $\F_p^4$ and let $B = [p/10]^4 \subset \F_p^4$ be a standard box. For a random shift $s$, we have $|(S' - s) \cap (B\times B)| \gtrsim q^{5/2}$. Let $\psi: B \rightarrow [n]$ be the map $\psi: (b_0, b_1, b_2, b_3) = \sum p^i b_i$ and define
    $$
    S = \psi((S' - s) \cap (B\times B)) \subset [n]^2.
    $$
    It is clear that $|S| \sim q^{5/2} \sim n^{5/4}$ and that $\psi$ preserves the property of being skew corner-free. So $S$ is skew corner-free.
    
\section{Proof of Theorem 1.2}

Let $X = \pi(S) \subset \mathbb{F}_{q}^{n}$ denote the projection of $S$ onto the first coordinate. For each $x \in X$, let $C_{x} \subset \mathbb{F}_{q}^{n}$ denote the set of elements $y \in \mathbb{F}_{q}^{n}$ with the property that $(x,y)$ is a point in $S$. Morally, the elements of $C_{x}$ identify the points of $S$ in the `column' above $x$, so we will sometimes refer to the set $C_{x}$ as the {\it{column above $x$}}. Clearly, $|S| = \sum_{x \in X} |C_{x}|$. In order to bound this sum we are going to use the so-called Croot-Lev-Pach lemma, in the same style as Ellenberg and Gijswijt used it in \cite{eg}. We first briefly recall the statement and introduce some useful notation for later.

Let $\mathcal{V}(q,n)$ be the $\mathbb{F}_{q}$-vector space of functions $f :\ \mathbb{F}_{q}^{n} \to \mathbb{F}_{q}$. A basis for this vector space is given by the set of monomials
$$M(q,n) = \left\{x_{1}^{a_1} \ldots x_{n}^{a_{n}}:\ 0 \leq a_i \leq q-1\right\}.$$
Given a positive integer $d$, let $M_{d}(q,n)$ be the set of monomials in $M(q,n)$ of degree at most $d$, and let $\mathcal{V}_{d}(q,n) \subset \mathcal{V}(q,n)$ be the set of polynomials of degree at most $d$ over $\mathbb{F}_{q}$ spanned by these monomials. Finally, let $m_{d}(q,n) = |M_{d}(q,n)|$. Using this terminology, the general form of the Croot-Lev-Pach lemma over $\mathbb{F}_{q}^{n}$ can be stated as follows.

\begin{lemma}\label{CLP}
Let $f \in \mathcal{V}_{d}(q,n)$ and let $A$ denote the $q^{n} \times q^{n}$ matrix with entries $A_{y,z} = f(y,z)$ for $y,z \in \mathbb{F}_{q}^{n}$. Then, $\operatorname{rank}(A) \leq 2 \cdot m_{d/2}(q,n)$. 
\end{lemma}

See for example \cite{pp}. Returning to skew corner-free sets, for a given positive integer $d$ whose value we will decide upon later, we now let $V \subset  \mathcal{V}_{d}(q,n)$ be the $\mathbb{F}_{q}$-space of polynomials vanishing on the complement of $X$. Note that the dimension of $V$ satisfies $\operatorname{dim}(V) \geq m_{d}(q,n) - q^{n} + |X|$. 

Let $P \in V$ be an element with the support $\Sigma:= \left\{x \in \mathbb{F}_{q}^{n}: P(x) \neq 0\right\}$ of maximum size. Note that $|\Sigma| \geq \operatorname{dim}(V)$ holds, since otherwise there would exist a nonzero $Q \in V$ vanishing on $\Sigma$. Such a polynomial $Q$ would generate an element of $V$ with larger support than $P$: indeed, notice that $(P+Q)(x) = P(x) \neq 0$ holds for every $x \in \Sigma$ and $(P+Q)(x) = Q(x) \neq 0$ must hold for some $x \not\in \Sigma$. 

Now, let $P \in V$ be a polynomial with support $\Sigma$ satisfying \begin{equation} \label{support}
|\Sigma| \geq \operatorname{dim}(V) \geq m_{d}(q,n) - q^{n} + |X|.
\end{equation}
For every element $x \in  X$, note that there are no distinct elements $y,z \in C_{x}$ such that $x+z-y \in  X$. Indeed, notice that this would yield a skew corner of the form $(x, y)$, $(x, y+d)$, $(x+d, y')$, where $d=z-y$ and $y'$ is some element in the column $C_{x+z-y}$ (which is non-empty if $x+z-y \in X$). 
In particular, if $x \in X$ is such that $P(x) \neq 0$, then the $q^{n} \times q^{n}$ matrix $A$ with rows and columns indexed by the elements of $\mathbb{F}_{q}^{n}$ and with entries $A_{y,z}=P(x+y-z)$ for $y,z \in \mathbb{F}_{q}^{n}$ has a very nice property: its the restriction to the set rows and columns corresponding to the elements of the column $C_{x}$ is a diagonal matrix with non-zero entries on the diagonal. By Lemma \ref{CLP}, it thus follows that
\begin{equation} \label{C_x} 
|C_{x}| \leq \operatorname{rank}(A) \leq 2m_{d/2}(q,n)
\end{equation}
holds for every $x \in X$ such that $P(x) \neq 0$. On the other hand, by \eqref{support}, the number of elements $x \in X$ with $P(x) = 0$ is
$$
|X| - |\Sigma| \leq q^{n} - m_{d}(q,n).
$$
Like in \cite{eg}, we next note that the quantity $q^{n} - m_{d}(q,n)$ represents the number of $q$-power-free monomials whose degree is greater than $d$, and these are in a simple bijection with the set of monomials of degree less than $(q-1)n-d$. Thus, $q^{n} - m_{d}(q,n)=m_{(q-1)n-d}(q,n)$. For each of these $x \in X$ with $P(x) = 0$, we shall use the trivial bound $|C_{x}| \leq q^{n}$. Putting things together and using the fact that $|\Sigma| \leq q^{n}$, we get that
$$\sum_{x \in X} |C_{x}| \leq |\Sigma| \cdot 2m_{d/2}(q,n) + m_{(q-1)n-d}(q,n) \cdot q^{n} \leq q^{n}\left(2m_{d/2}(q,n) + m_{(q-1)n-d}\right).$$
Picking $d = 2(q-1)n/3$, we get that
$$\sum_{x \in X} |C_{x}| \leq q^{n} \cdot 3m_{(q-1)n/3}(q,n).$$
Since
$$m_{(q-1)n/3}(q,n) \leq \inf_{0< x < 1} x^{-(q-1)n/3}(1+x+\ldots+x^{q-1})=q^{(1-c_{q})n}$$
(see for example \cite{eg}), the conclusion follows.

\bibliographystyle{amsplain0.bst}
\bibliography{main}

\end{document}